\documentclass{amsart}

\usepackage{xspace,amssymb,amsfonts,euscript,eufrak,mathrsfs}

\usepackage{amsthm,amsmath}
\usepackage{url}
\usepackage{graphicx}
 \usepackage{epstopdf}
\usepackage{xcolor}
\usepackage{color}
\usepackage[all]{xy}

\newcommand{\BMP}{\mathscr B}

\newcommand{\E}{\mathcal{E}}
\newcommand{\F}{\mathcal{F}}
\newcommand{\h}{\mathfrak{h}}
\newcommand{\I}{\mathcal{I}}

\newcommand{\MG}{\mathcal{G}}

\newcommand{\MGJ}{\mathcal{G}^{J}}
\newcommand{\rk}{\underline{\text{rk}}}
\newcommand{\SR}{\mathcal{S}}
\newcommand{\T}{\mathcal{T}}
\newcommand{\V}{\mathcal{V}}
\newcommand{\W}{\mathcal{W}}
\newcommand{\Z}{\mathcal{Z}}

\newcommand{\g}{\mathfrak{g}}

\renewcommand{\check}[1]{{#1}^{\vee}}

\newtheorem{prop}{Proposition}[section]
\newtheorem{theor}{Theorem}[section]
\newtheorem{defin}{Definition}[section]
\newtheorem{cor}{Corollary}[section]
\newtheorem{lem}{Lemma}[section]

\newtheorem{quest}{Question}[section]
\newtheorem{rem}{Remark}[section]

\begin{document}
\title{Kazhdan-Lusztig combinatorics in the moment graph setting}
\author{Martina Lanini}
\address{
Department Mathematik, Friedrich-Alexander-Universit\"{a}t Erlangen-N\"{u}rnberg, Cauerstr. 11,  91058  Erlangen, Germany}
\email{lanini@mi.uni-erlangen.de}
\maketitle
\vspace{-9mm}
\begin{abstract}
Motivated by a question on the graded rank of the stalks of the canonical sheaf on a Bruhat graph, we lift some equalities concerning (parabolic) Kazhdan-Lusztig polynomials to this moment graph setting. Our proofs hold also in positive characteristic,  under some technical assumptions.
\end{abstract}

\section{Introduction}
In 1979 Kazhdan and Lusztig (\cite{KL}) associated to a given Coxeter group $\W$ a family of polynomials $\{P_{y,w}(q)\}$ indexed by pairs of elements in $\W$. In the case $\W$ was a Weyl group, then $P_{y,w}(q)$ was related 
to the local intersection cohomology of the corresponding Schubert variety(cf.Appendix A of \cite{KL} and \cite{KL80}). Some years later,  Deodhar in \cite{Deo87} introduced the parabolic analogue of Kazhdan-Lusztig polynomials. 
Namely, if $(\W, \SR)$ is a Coxeter system, $J\subseteq \SR$, and $\W^J$ is the set of minimal coset representatives of  $\W/\langle J \rangle$, he defined two families of polynomials $\{P^{J,-1}_{y,w}(q)\}$ and  $\{P^{J,q}_{y,w}(q)\}$, where $y,w\in\W^J$. The $\{P^{J,-1}_{y,w}(q)\}$ are a generalisation of the polynomials defined by Kazhdan and Lusztig and for $J=\emptyset$ they coincide. As in the regular case, if $\W$ is a Weyl group, then these polynomials have a geometrical meaning, and, in particular, they are related to the intersection cohomology  of the corresponding partial Schubert variety.

Kazhdan and Lusztig  (\cite{KL}), resp. Lusztig (\cite{Lu80c}), conjectured that the Kazhdan-Lusztig polynomials played a very important role in the representation theory of complex Lie algebras, resp. of semisimple, 
simply connected, reductive algebraic groups over a field of positive characteristic. The characteristic zero setting is now well understood (cf.\cite{KL80},\cite{BeBe},\cite{BK}), while the positive characteristic analogue is not. Actually Lusztig's conjecture was \emph{almost} proved in the 90s  via the joint work of Kazhdan-Lusztig (\cite{KL93}), Kashiwara-Tanisaki (\cite{KT}) and Andersen-Jantzen-Soergel (\cite{AJS}). Here \emph{almost} means that it was possible to prove the conjecture only  if the characteristic of the base field is big enough, since it was obtained as a limit of the characteristic zero case. A new approach to Lusztig's conjecture is due to Fiebig (\cite{Fie08a},\cite{Fie07a}) and it is based on the theory of sheaves on moment graphs. 

Moment graphs were introduced by Goresky, Kottwitz and MacPherson (\cite{GKM}), in order to study the equivariant cohomology of a complex algebraic variety equipped with a torus action and having some nice properties. 
 In 2001 Braden and MacPherson in \cite{BM01} were able to describe the equivariant intersection cohomology of such a variety via sheaves on the moment graph. In particular, if $\W$ is a Weyl group with $\SR$, the set of simple reflections, and $J\subseteq \SR$, Braden and MacPherson associated to $w\in \W^J$ a sheaf $\BMP^J_w$: the \emph{canonical} --- or \emph{BMP} --- sheaf. 
This  object describes the local  intersection cohomology of the corresponding Schubert variety in a partial flag variety. Braden-MacPherson's construction was performed in characteristic zero, but it is possible 
to develop this theory in any characteristic. Fiebig and Williamson proved in \cite{FW} that, with certain technical assumptions, 
in positive characteristic $\BMP_w^J$ computes the stalks of indecomposable parity sheaves (introduced in \cite{JMW}). 
It is  now natural to ask whether it is possible to connect the canonical sheaf to  Kazhdan-Lusztig polynomials.

\vspace{2mm}
\textbf{Question \ref{mult_conj}.}(cf.\cite{Fie07b}, Conjecture 4.4) \emph{Under which assumptions on the characteristic of the base field, do we have   $\rk \,(\BMP_w^J)^y=P^{J,-1}_{y,w}$ for $y\leq w$  and $y,w$ varying in some relevant subset of $\W^J$?}
\vspace{2mm}

This equality is true in characteristic zero for any pair $y,w$ and in this case it is equivalent to  Kazhdan-Lusztig's conjecture (cf.\cite{Fie08a}). In characteristic $p$, if $\W$ is affine and if we only consider $w$ restricted (cf.\cite{Fie07a}), it is proved for $p$
bigger than a huge (but  explicit) bound (cf.\cite{Fie08c}), and, for $p$ bigger than the Coxeter number,
 it implies Lusztig's conjecture (cf.\cite{Fie07b}, \cite{Fie07a}). From a recent result of Polo (private communication, 7 May, 2012), it follows that 
if $\W=S_{4p}$ the stalks of the BMP-sheaf are definitively not given by these polynomials (see \S\ref{ssec_mult_conj} for more details). Anyway, this question motivates our work, since it now makes 
sense to interpret some equalities concerning  (parabolic) Kazhdan-Lusztig polynomials in terms of stalks of the canonical sheaves. In order to lift properties of KL-polynomials to the level of canonical sheaves, we will use two different techniques: 
the pullback of canonical sheaves (see Section \ref{section_pullbacks}) and an action of the Weyl group on the set of global sections of the BMP-sheaf (see Section \ref{section_invariants}).

Let $k$ be a local ring with $2\in k^{\times}$. We define the notion of $k$-homomorphism between two moment graphs and of pullback of sheaves. These will provide a useful tool, namely, under some assumptions on $k$,

\vspace{2mm}
\textbf{Lemma \ref{pullbackBMP}.} \emph{Let $\MG$ and $\MG'$ be two moment graphs, both with a unique maximal vertex, w resp. w', and let $f$ be a $k$-isomorphism between them. If  $\BMP_w$ and $\BMP'_{w'}$ are the corresponding canonical sheaves, then $\BMP_w\cong f^*\BMP'_{w'}$ as $k$-sheaves on $\MG$.}
\vspace{2mm}

Thanks to this result, in some good situations it will be enough to study the combinatorics of the underlying moment graphs that in our case are just labeled, oriented Bruhat graphs (see \S\ref{bruhatMG}). This is the case in the following theorem:

\vspace{2mm}
\textbf{Theorem \ref{KLpropsthm}.} \emph{Let $y,w\in \W$ be such that $y\leq w$, then 
\begin{itemize}
\item[(i)] $\BMP_w^y\cong \BMP_{w^{-1}}^{y^{-1}}$. 
\end{itemize}
Let $s\in \SR$ be such that $ws<w$, then 
\begin{itemize}
\item[(ii)] if $y\not\leq ws$, $\BMP_w^y\cong\BMP_{ws}^{ys}$,
\end{itemize}}\vspace{2mm}
where we write $\BMP_w$ instead of $\BMP_w^{\emptyset}$.

 The last part of the paper is devoted to the study of an action of  a certain subgroup of the Weyl group $\W$ on the space of global sections of the canonical sheaf and, in particular, to the proof  that the data we need to  build the canonical sheaf is contained in the invariants with respect to this action. This result, together with some combinatorics of the corresponding Bruhat graph, gives us the categorical analogue of a result due to Kazhdan and Lusztig (cf.\cite{KL80}):

\vspace{2mm}
\textbf{Theorem \ref{Bys}.} \emph{Under some assumptions on $k$, if $y,w\in \W$ and $s\in \SR$ are such that  $y\leq w$ and $ws<w$, then 
$\BMP_w^y\cong \BMP_w^{ys}$.}
\vspace{2mm}

Inspired by a theorem of Deodhar (\cite{Deo87}) we prove a relation between the canonical sheaf on a regular Bruhat graph $\MG$ and the ones on the corresponding  parabolic Bruhat graphs $\MGJ$, for $J$ such that the subgroup $\langle J \rangle$ is finite. Let $w_J$ be the longest element  of $\langle J \rangle$ then, under some assumptions on $k$, we have

\vspace{2mm}
\textbf{Theorem \ref{BMPparab}.} \emph{If $y,w\in \W^J$, then $(\BMP_w^J)^y\cong (\BMP_{ww_J}^{\emptyset})^{yw_J}$.}
\vspace{2mm}

In order to prove this result we consider again  the space of invariants with respect to the above action. The claim follows from the connection between  this module and the parabolic canonical sheaf.

\subsection*{Structure of the paper. } Sections 2 and 3 are about moment graphs and sheaves on them. In Section 4 we introduce Braden-MacPherson sheaves and recall some of their properties. We develop and apply  the technique of pullbacks in Section 5, while the one of invariants is used in the last section.

\section{Moment graphs}

In this section we recall the definition of moment graphs on a lattice and we define the notion of $k$-homomorphism between two moment graphs. 

Let $k$ be from now on a local ring inside which $2$ is an invertible element. 
 Let  $Y\cong \mathbb{Z}^r$ be a lattice of finite rank and denote by $Y_k:=Y\otimes_{\mathbb{Z}} k$ .

\begin{defin}A \emph{moment graph} $\MG$ on $Y$ is given by $(\V,\E, \trianglelefteq, l)$, where:
\item[(i)] $(\V,\E)$ is a directed graph without directed cycles nor multiple edges,
\item[(ii)] $ \trianglelefteq$ is a partial order on $\V$ such that if $x,y\in \V$ and $E:x\rightarrow y\,\in \E$ then $x \trianglelefteq y$,
\item[(iii)] $l:\E\rightarrow Y\setminus \!\{0\}\,$  is a map  called the \emph{label function}.
\end{defin}

\begin{defin}\label{def_GKM}Let $\MG$ be  a moment graph on the lattice $Y$, then 
 \begin{itemize}
\item $\MG$ is a \emph{$k$-moment graph} on $Y$ if all labels are non-zero in $Y_k$
\item $(\MG, k)$ is a\emph{ $GKM$-pair} if all pairs of distinct edges containing a common vertex have labels $k$-linearly independent in $Y_k$. 
\end{itemize}
\end{defin}

Observe that if $(\MG,k)$ is a $GKM$-pair, then $\MG$ is a $k$-moment graph. These  properties are  very important and in the sequel  they  will give  a restriction on the ring $k$.

\subsection{$k$-homomorphisms of moment graphs.}\label{homoMG} Let $\MG=(\V,\E,  \trianglelefteq, l)$ and  $\MG'=(\V',\E', \trianglelefteq', l')$ be two moment graphs on  $Y$. Since a moment graph is given by an oriented and ordered graph plus some other data coming from $Y$, we define a $k$-homomorphism as a map of graphs plus a collection of automorphisms of the $k$-module $Y_k$ satisfying certain requirements. More precisely,

\begin{defin}\label{Def_morphMG}A $k$-\emph{homomorphism} between two moment graphs on $Y$
\begin{equation*}f:(\V,\E, \trianglelefteq, l)\rightarrow (\V',\E', \trianglelefteq', l')
\end{equation*} is given by $\left(f_{\V}, (f_{l,x})_{x\in \V}\right)$, where
\item[(MORPH1)] $f_{\V}:\V \rightarrow \V'$ is any (order preserving) map of posets such that, 
if $x-\!\!\!-\!\!\!-y\in\E$, then either $f_{\V}(x)-\!\!\!-\!\!\!-f_{\V}(y)\in \E' $, or $f_{\V}(x)=f_{\V}(y)$.
\item[(MORPH2)] For all $x\in \V$, $f_{l,x}:Y_k\rightarrow Y_k\in\text{Aut}_k(Y_k)$  is such that,
 if $E:x-\!\!\!-\!\!\!-y\in \E$ and $f_{\V}(x)\neq f_{\V}(y)$, the following two conditions are verified:
\hspace{5mm}\item[(MORPH2a)] $f_{l,x}(l(E))=h\cdot l'(f_{\E}(E)) $, for some $h\in k^{\times}$
\hspace{5mm}\item[(MORPH2b)]  $\pi \circ f_{l,x} =\pi \circ f_{l,y}$,
where $\pi$ is the canonical quotient map $\pi:Y_k\rightarrow Y_k/l'(f_{\E}(E))Y_k$.
\end{defin}

For an edge $E:x-\!\!\!-\!\!\!-y\in \E$ such that $f_{\V}(x)\neq f_{\V}(y)$, we will denote $f_{\E}(E):= f_{\V}(x)-\!\!\!-\!\!\!-f_{\V}(y)$.

\begin{defin} $f=\left(f_{\V},(f_{l,x})_{x\in\V}\right):\MG\rightarrow\MG'$ is a \emph{$k$-isomorphism} of moment graphs if  $f:\MG\rightarrow\MG'$ is a $k$-homomorphism and  the following two conditions hold:
\item[(ISO1)] $f_{\V}$ is an isomorphism of posets
\item[(ISO2)] for all $u\rightarrow w\in \E'$, there exists exactly one $x\rightarrow y\in \E$ such that $f_{\V}(x)=u$ and $ f_{\V}(y)=w $. 
\end{defin}

If $f$ is a $k$-isomorphism from the moment graph $\MG$ to itself, we say that it is a \emph{$k$-automorphism} of $\MG$.

\subsection{Bruhat graphs. }\label{bruhatMG} Here we describe a class of moment graphs, that is for our purposes the most important one. 

 We start by recalling some notation from \cite{Kac}.  Let $\g$ be a symmetrisable Kac-Moody algebra, that is the Lie algebra $\g (A)$ associated to a symmetrisable generalised Cartan matrix $A$, and $\h$ its Cartan subalgebra. Let $\Pi=\{\alpha_i\}_{i=1,\ldots , n}\subset \h^*$, resp. $\check{\Pi}=\{\check{\alpha_i}\}_{i=1,\ldots , n}\subset \h$, be the set of simple roots, resp. coroots; let  $\Delta$, resp. $\Delta_+$, resp. $\Delta_{+}^{\text{re}}$  be the root system, resp. the set of positive roots, resp. the set of positive real roots; and let $Q=\sum_{i=1}^{n}\mathbb{Z}\alpha_i$, resp. $\check{Q}=\sum_{i=1}^{n}\mathbb{Z}\check{\alpha_i}$, be the root lattice, resp. the coroot lattice. For any $\alpha\in\Delta$, we denote by $s_{\alpha}\in GL(\h^*)$ the reflection, whose action on $v\in\h^*$ is given by
\begin{equation} s_{\alpha}(v)=v-\langle v, \check{\alpha} \rangle \alpha
\end{equation}
 Let $\W=\W(A)$ be the Weyl group associated to $A$, that is the subgroup of $GL(\h^*)$ generated by the set of simple  reflections $\SR=\{s_{\alpha}\mid \alpha\in\Pi\}$. Recall that $(\W, \SR)$ is a Coxeter system (cf.\cite{Kac}, \S 3.10).

However, $\W$ can be seen also as a subgroup of $GL(\h)$, by setting, for any $\lambda\in \h$
\begin{equation}
s_{\alpha}(\lambda):=\lambda-\langle \alpha ,  \lambda \rangle \check{\alpha}
\end{equation}

We will denote by $\T\subset \W$ the set of reflections, that is 
\begin{equation}\T=\left\{  s_{\alpha}\mid \alpha\in\Delta_{+}^{\text{re}} \right\}=\left\{wsw^{-1} \mid w\in\W, \,\,s\in\SR \right\}
\end{equation}
Hereafter we will write $\alpha_t$ to denote the positive real root corresponding to the reflection $t\in\T$.
Finally, denote by $\ell:\W\rightarrow \mathbb{Z}_{\geq 0}$ the length function and by $\leq$ the Bruhat order on $\W$. 

For any $J\subseteq \SR$, denote by $\W_{J}:=\langle J \rangle$ and by $\W^J$ the set of minimal coset representatives of the equivalence classes of $\W/\W_J$. 

\begin{defin}\label{def_parBruhatMG}Let $\W$, $\SR$ and $J$ be as above. Then the \emph{Bruhat (moment) graph} $\MGJ=\MG(\W^J)=(\V, \E, \leq, l)$ associated to $\W^J$ is a moment graph on $\check{Q}$ and it is given by
\item[(i)] $\V=\W^J$
\item[(ii)] $\E=\left\{ x\rightarrow y \mid x< y\,,\, \exists\, \alpha\!\in\!\Delta_+^{\text{re}}, \, \exists  w\in\W_J \text{ such that } ywx^{-1}=s_{\alpha} \right\}$
\item[(iii)] $l(x\rightarrow s_{\alpha}xw^{-1}):=\check{\alpha}$.
\end{defin}

Such a moment graph has an important geometric meaning. Indeed, there is a partial flag variety $Y$ corresponding to $\W$ and $J$ as above (see \cite{Ku}) and it carries an action of a torus $T$, whose Lie algebra is $\h$, and a ($T$-invariant) stratification with certain good properties (see \cite{BM01}). The Bruhat graph encodes the action of this torus, in particular, the vertices are the $0$-dimensional orbits, while  the edges represent the $1$-dimensional orbits (cf.\S 2.1 of \cite{FW}). The  partial order on the set of vertices is induced by  the  stratification coming from  the  decomposition $Y=\bigsqcup_{w\in \W^J} Y_w$, where, indeed,  $\overline{Y_w}=\bigsqcup_{\substack{y\leq w\\y\in \W^J}} Y_y$.

\section{Category of $k$-sheaves on a moment graph}

Consider a moment graph $\MG=(\V,\E,\trianglelefteq,l)$ on a lattice $Y$. Recall that for any local  ring $k$  with  $2\in k^{\times}$ we denoted by $Y_k:=Y\otimes_{\mathbb{Z}}k$.

Let $S_k:=\text{Sym}(Y_k)$ be the symmetric algebra of $Y_k$. We provide $S_k$ with a  $\mathbb{Z}$-grading such that $(S_k)_{\{2\}}=Y_k$. From now on every $S_k$-module will be finitely generated and $\mathbb{Z}$-graded and every morphism between $S_k$-modules will be of degree zero.

\begin{defin}A \emph{$k$-sheaf} $\mathcal{F}$ on  $\MG$ is given by $(\{\mathcal{F}^{x}\},\{\mathcal{F}^E\}, \{\rho_{x,E}\})$, where:
\item[(i)   ] for all $x\in\V$, $\mathcal{F}^{x}$ is an $S_k$-module;
\item[(ii) ] for all $E\in\E$, $\mathcal{F}^{E}$ is an $S_k$-module such that $l(E)\cdot \mathcal{F}^E=\{0\}$;
\item[(iii)] for $x\in \V$ and $E\in \E$ such that $x$ is in the boundary of the edge $E$, the map $\rho_{x,E}:\mathcal{F}^x\rightarrow \mathcal{F}^E$ is  a homomorphism of $S_k$-modules. 
\end{defin}

\begin{defin} A \emph{homomorphism} $\varphi:\mathcal{F}\longrightarrow \mathcal{F'}$ between $k$-sheaves  on the moment graph $\MG$ is given by the following data
\item[(i)   ] for all $x\in\V$, $\varphi^x:\mathcal{F}^{x}\rightarrow \mathcal{F'}^x$ is a homomorphism of $S_k$-modules
\item[(ii) ] for all $E\in\E$,  $\varphi^E:\mathcal{F}^{E}\rightarrow \mathcal{F'}^E$ is a homomorphism of $S_k$-modules
such that  for any $x\in\V$ on the border of $E\in\E$ the following diagram commutes
\begin{equation*}
\begin{xymatrix}{
\mathcal{F}^x \ar[d]_{\rho_{x,E}}\ar[r]^{\varphi^x}&\mathcal{F'}^x \ar[d]^{\rho'_{x,E}}\\
\mathcal{F}^E \ar[r]_{\varphi^E}&\mathcal{F'}^E\\
}
\end{xymatrix}
\end{equation*}
\end{defin}

We denote by  $\frak{Sh}_{\MG}^k$  the category of \emph{$k$-sheaves} on $\MG$ having as objects the $k$-sheaves on $\MG$ and as morphisms the homomorphisms between them.

\subsection{Pullback sheaves.} Let $\MG=(\V, \E, \trianglelefteq, l)$ end $ \MG'=(\V',\E',\trianglelefteq', l')$ be two moment graphs on $Y$ and fix
 $f:\MG\longrightarrow \MG'$ a $k$-homomorphism of moment graphs on $Y$(cf.\S \ref{homoMG}).

\begin{defin}\label{Def_pull} Let $\F\in Ob(\frak{Sh}_{\MG'}^k)$, then $f^*\F\in  Ob(\frak{Sh}_{\MG}^k)$ is defined as follows
\item[(PULL1)] for all $x\in\V$, $(f^*\F)^x:=\F^{f_{\V}(x)}$ and $s\in S_k$ acts on it via $f_{l,x}(s)$
\item[(PULL2)] for all $E:x-\!\!\!-\!\!\!-y\in\E$\begin{equation*}(f^*\F)^{E}=\left\{ \begin{array}{ll}
\F^{f_{\V}(x)}/l(E) \F^{f_{\V}(x)} & \text{if } f_{\V}(x)=f_{\V}(y)\\
\F^{f_{\E}(E)} &\text{ otherwise }
\end{array}
\right.
\end{equation*}
and of $s\in S_k$ acts on $(f^*\F)^{E}$ via $f_{l,x}(s)$.
\item[(PULL3)] for all $x\in \V$ and $E\in \E'$, such that $E:x-\!\!\!-\!\!\!-y$,\begin{equation*}(f^*\rho)_{x,E}=\left\{ \begin{array}{ll}
\text{canonical quotient map} & \text{if } f_{\V}(x)=f_{\V}(y)\\
\rho_{f_{\V}(x),f_{\E}(E)} &\text{ otherwise }
\end{array}
\right.
\end{equation*}
\end{defin}

\begin{rem}For all $E\in\E$, the action of $S_k$ on $(f^*\F)^{E}$ in \textit{(PULL2)} is well-defined thanks to conditions \textit{(MORPH2a)} and \textit{(MORPH2b)}.
\end{rem}

We say that $f^*\mathcal{F}$ is the \emph{pullback} of $\mathcal{F}$. In what follows, the notion of pullback sheaf will allow us to compare $k$-sheaves on different moment graphs.

\subsection{Sections of sheaves.}\label{sections_sheaves}  For each $\mathcal{I}\subset \V$ we can consider the \emph{set of  local sections} of a $k$-sheaf $\mathcal{F}\in Ob(\mathfrak{Sh}_{\MG}^k)$ over $\mathcal{I}$:
\begin{equation*}\Gamma(\mathcal{I},\mathcal{F}):=\left\{ (m_x)\in \prod_{x\in \mathcal{I}} \mathcal{F}^{x} \mid \begin{array}{c} 
\rho_{x,E}(m_x)=\rho_{y,E}(m_y)\\ \forall E:x-\!\!\!-\!\!\!- y\in\E ,\, x,y\in\mathcal{I}                                                                                                                  
                                                                                                                 \end{array}
   \right\}.
\end{equation*}
We denote by $\Gamma(\mathcal{F})=\Gamma_k(\V,\mathcal{F})$ the set of global sections of the $k$-sheaf $\mathcal{F}$.

We call \emph{$k$-structure algebra} of the moment graph $\MG$ the set
\begin{equation*}\Z_k=\Z_k(\MG):=\left\{ (z_x)\in \prod_{x\in \V}S_k \mid  \begin{array}{c} z_x-z_y\in l(E)S_k  \\ \forall  E:x-\!\!\!-\!\!\!- y\in\E                                                                                                                  
                                                                                                                \end{array}  \right\}.
\end{equation*}
It is easy to check that for any $\mathcal{F}\in Ob(\mathfrak{Sh}_{\MG}^k)$ the $k$-structure algebra $\Z_k$ acts on $\Gamma(\mathcal{F})$ via componentwise multiplication. 

\subsection{Flabby sheaves. } We use the order on the set of vertices of a moment graph $\MG$ to define a topology on it: the \textit{Alexandrov topology}. We say that $\mathcal{I}$ is open if for any $x\in \mathcal{I}$ and any $y\in \V$ such that $x\trianglelefteq y$ then $y\in  \mathcal{I}$ as well. 

A classical question in sheaf theory is to ask whether a sheaf is flabby or not, that is whether any local section over an open set extends to a   global one or not.

Let $\mathcal{F}\in Ob(\mathfrak{Sh}_{\MG}^k)$. We fix a vertex $x\in \V$ and we denote
\begin{equation*}\E_{\delta x}:=\left\{  E\in\E\mid E:x\rightarrow y         \right\}
\end{equation*}
 \begin{equation*}\V_{\delta x}:=\left\{ y\in \mathcal{V}\mid \exists E\in\E_{\delta x}  \text{ such that } E:x\rightarrow y       \right\}.
   \end{equation*}

Now we define $\mathcal{F}^{\delta x}$ as the image of $\Gamma(\{ \triangleright x \}, \mathcal{F} ):=\Gamma(\{ y\in \V\,\,|\,\, y\triangleright x   \},\mathcal{F})$ under the composition of the following functions:
\begin{equation*}\begin{xymatrix}{
u_x:\Gamma(\{\triangleright x\},\mathcal{F})\ar[r]&\bigoplus_{y\triangleright  x}\!\mathcal{F}^y\ar[r]&\bigoplus_{y\in \V_{\delta x}}\!\!\mathcal{F}^y\ar[r]^{\oplus \rho_{y,E}}&\bigoplus_{E\in \E_{\delta x}} \!\!\mathcal{F}^E\\
}\end{xymatrix}
\end{equation*}

Denote 
\begin{equation*}\begin{xymatrix}{
d_x=\oplus_{E\in\E_{\delta x}} \rho_{x,E}:\F^x\ar[r]&\bigoplus_{E\in \E_{\delta x}} \!\!\mathcal{F}^E\\
}\end{xymatrix}
\end{equation*}

Observe that $m\in\Gamma(\{\triangleright x\}, \F)$ can be extended, via $m_x$, to a section $\tilde{m}=(m,m_x)\in\Gamma(\{\trianglerighteq x\}, \F)$ if and only if $d_x(m_x)=u_x(m)$. This fact motivates the following result, due to Fiebig, that gives a characterization of the flabby objects in $\frak{Sh}_{\MG}^k$.

\begin{prop}[\cite{Fie08a}]\label{prop_flabby} Let $\mathcal{F}\in Ob(\mathfrak{Sh}_{\MG}^k)$. Then the following are equivalent:
\begin{itemize}
\item[(i)] $\mathcal{F}$ is flabby with respect to the Alexandrov topology, i.e.  for any open $\mathcal{I}\subseteq \V$ the restriction map $\Gamma(\mathcal{F})\rightarrow \Gamma(\mathcal{I}, \mathcal{F})$ is surjective.
\item[(ii)] For  any vertex $x\in \V$ the restriction map $\Gamma(\{\trianglerighteq x\},\mathcal{F})\rightarrow \Gamma(\{\triangleright x\}, \mathcal{F})$ is surjective.
\item[(iii)] For any vertex $x\in \V$  the map $d_x:\mathcal{F}^x\rightarrow \bigoplus_{E\in \E_{\delta x}}\mathcal{F}^E$ contains  $\mathcal{F}^{\delta x}$ in its image.
\end{itemize}
\end{prop}

\section{Braden-MacPherson sheaves.} 

In this section we introduce the most important object of our paper, namely the canonical sheaf. It was first defined by Braden and MacPherson --- only in characteristic zero --- in order to compute  certain  intersection cohomology complexes. Despite  this, their algorithm  makes sense in any characteristic and Fiebig and Williamson proved   in \cite{FW} that it computes the multiplicities of parity sheaves (see \cite{JMW}) in positive characteristic if $(\MG, k)$ is a \emph{GKM}-pair. The following theorem allows us to consider this  sheaf.

\begin{theor}[\cite{BM01}, char $k$=0; \cite{Fie08a}] \label{BM01}Let $\MG$ be a finite $k$-moment graph over $Y$ with highest vertex $w$. There exists exactly one (up to isomorphism) indecomposable $k$-sheaf $\BMP_w$ on $\MG$ with the following properties:  
\begin{itemize}
\item[(i)] $\BMP_w^w\cong S_k$;
\item[(ii)] If $x,y\in \V$, $E:x\rightarrow y\in \E$, then the map  $\rho_{y,E}:\BMP_w^y\rightarrow \BMP_w^E$ is surjective with kernel  $l(E)\BMP_w^y$;
\item[(iii)]  If $x,y\in \V$, $E:x\rightarrow y\in E$, then $\rho_{\delta x}:=\bigoplus_{E\in \E_{\delta x}}\rho_{x,E}: \BMP_w^x\rightarrow \BMP_w^{\delta x}$ is a projective cover in the category of graded $S_k$-modules.
\end{itemize}
\end{theor}

We call  $\BMP_w$  the \textit{Braden-MacPherson sheaf} or the \textit{canonical sheaf}. We will also refer to it as the \emph{BMP-sheaf}.

\begin{rem} By Theorem \ref{BM01} and Proposition \ref{prop_flabby}, the canonical sheaf is flabby for the Alexandrov topology. This property will  be crucial in what follows.
\end{rem}

\subsection{Graded rank of the stalks of a BMP-sheaf.}\label{ssec_mult_conj} For $j\in \mathbb{Z}$ and $M$ a graded $S$-module we denote by $M\{j\}$ the graded $S$-module obtained from $M$ by shifting the grading by $j$, i.e. $M\{j\}_{\{i\}}=M_{\{j+i\}}$. If $M=\bigoplus_{i=1}^n S_k\{j_i\}$, then its graded rank is $\rk\, M=\sum_{i=1}^n q^{-\frac{j_i}{2}}\in \mathbb{Z}_{\geq 0}[q^{\frac{1}{2}},q^{-\frac{1}{2}}]$.

 Let $\MGJ$ be the Bruhat graph we defined  in \S\ref{bruhatMG}. Thus for any $w\in \W^J$ we can consider the subgraph $\MG^J_{w}:=\MGJ_{|\{\leq w\}}$. It is a finite $k$-moment graph (for any $k$) with highest vertex $w$, hence we may build the corresponding Braden-MacPherson sheaf $\BMP^J_w\in Ob(\mathfrak{Sh}_{\MG^J_w}^k)$ and we have:

\begin{quest}\label{mult_conj}Under which assumptions on the characteristic of the base field, do we have   $\rk \,(\BMP_w^J)^y=P^{J,-1}_{y,w}$ for $y\leq w$  and $y,w$ varying in some relevant subset of $\W^J$?
\end{quest}
Stand $P^{J,-1}_{y,w}$ for the parabolic Kazhdan-Lusztig polynomial (corresponding to the parameter $u=-1)$ introduced by Deodhar in \cite{Deo87}.

 If $k=\mathbb{Q}$, then $\rk \,(\BMP_w^J)^y=P^{J,-1}_{y,w}$  for any pair  $y, w\in \W^J$, with $y\leq w$, from \cite{KL}, \cite{KL80} and \cite{BM01}. Moreover, Fiebig  proved that in this case the equality is equivalent to a character formula conjectured by Kazhdan and Lusztig in \cite{KL} (see \cite{Fie08a}). Hence for a fixed pair of elements, from the characteristic zero case we get that the equality holds for $p$ large enough. Observe that the bound depends on the pair and there is no global bound in the infinite case.

Now let  $\W$ be an affine Weyl group, $h$ be its Coxeter number, $k$ be a field of characteristic $p\geq h$ and $y,w$ be restricted elements (cf.\cite{Fie07a}). A positive answer to Question \ref{mult_conj} would imply a conjecture by Lusztig (cf.\cite{Fie07b},\cite{Fie07a}). Fiebig was able to prove that the stalks $(\BMP_w^J)^y$  have the expected graded rank for $p$ bigger than an explicit --- but huge! --- bound depending on $\W$ (cf.\cite{Fie08c}). Motivated by the fact that in  the affine case the \emph{GKM}-condition (see Definition \ref{def_GKM}, (ii)) for the Bruhat graph of restricted elements is precisely $p\geq h$ ([12], Lemma 4.3), Fiebig suggested  the answer to Question \ref{mult_conj} to be  yes as soon as the \emph{GKM}-condition were satisfied (cf.\cite{Fie07b}, Conjecture 4.4).

Actually, very recently this conjecture has been  proven to be false for $\W=S_{4p}$. Indeed, Polo (private communication, 7 May, 2012) produced a family of elements $w_n$ in $S_{4n}$ (for each integer $n \geq 2$) such that there is $n$-torsion in some costalk of the intersection cohomology of the Schubert variety corresponding to $w_n$. The fact that this provides us with a family of counterexamples to Fiebig's conjecture is not immediate at all. We have to notice first that the Bruhat graph for $\mathfrak{sl}_r$ and $k$ constitute  a \emph{GKM}-pair for any $r$ and any field $k$ of characteristic $p>2$ and then to translate Question \ref{mult_conj} in terms of intersection cohomology complexes and parity sheaves (cf.\cite{FW}, Theorem 9.2).

Finally, let us consider an affine Weyl group  $\W$, whose Coxeter number is $h$ and a field $k$ of  characteristic $p>h$, but let us make $y,w$ vary in the finite Weyl group $\W^{\text{f}}<\W$. In this case a positive answer to Question \ref{mult_conj} would imply the Lusztig's conjecture around the Steinberg weight, which was presented by Soergel in the 90's as ``toy model'' for the original  Lusztig's conjecture (cf.\cite{WOb}).

In view of Polo's result, the bound $p\geq h$ seems to be the right one for expecting Question  \ref{mult_conj} to have a positive answer for $y,w$ restricted, resp, for any pair $y,w$, if $\W$ is affine, resp. finite. Moreover in this case, this problem  would still be related to Lusztig's conjecture, resp. Lusztig's conjecture on the Steinberg weight, as discussed above.

Anyway, Question \ref{mult_conj} proposes us an explicit formula connecting canonical sheaves and parabolic Kazhdan-Lusztig polynomials and motivates our work. We will indeed interpret in terms of stalks of BMP-sheaves some well-known identities concerning these polynomials.

\section{Pullback of BMP-sheaves.}\label{section_pullbacks}
The following lemma tells us that the pullback functor $f^*$ preserves canonical sheaves if $f$ is a $k$-isomorphism.
 \begin{lem} \label{pullbackBMP} Let $\MG$ and $\MG'$ be two $k$-moment graphs on $Y$, both with a unique maximal vertex, w resp. w', and let  $f:\MG\longrightarrow \MG'$ be a $k$-isomorphism. If  $\BMP_w$ and $\BMP'_{w'}$ are the corresponding canonical sheaves, then $\BMP_w\cong f^*\BMP'_{w'}$ as $k$-sheaves on $\MG$.
\end{lem}

\proof 
Let  $\MG=(\V,\E,\trianglelefteq,l)$, $\MG'=(\V',\E',\trianglelefteq',l')$ and $f=\left(f_{\V},(f_{l,x})\right)$. 

Notice that   $\I\subseteq \V$ is an open subset if and only if  $\I':=f_{\V}(\I)\subseteq \V  '$ is an open subset.  We prove that ${\BMP_w}_{|_{\I}}\cong {f^*\BMP'_{w'}}_{|_{\I'}}$ by induction on $|\I|=|\I'|$, for $\I$ open.

If  $|\I|=|\I'|=1$, we have $\I=\{w\}$ and  $\I'=\{w'\}$. In this case $\BMP_{w}^w=S_k$, $ \BMP_{w'}'^{w'}= S_k$ and the isomorphism $\varphi^w:\BMP_{w}^w\rightarrow \BMP_{w'}'^{w'}$ is just given by the twisting of the $S_k$-action, coming from  the automorphism of  $S_k$ induced by the automorphism  $f_{l,w}$  of $Y_k$.

Now let $|\I|=|\I'|=n>1$ and $y\in \I$ be a minimal element. Obviously, $y':=f_{\V}(y)$ is  also a minimal element for $\I'$. Moreover, for any $E\in\E$ we set   $E':=f_{\E}(E)$.

First of all, observe that $z \in \V_{\delta y}$ if and only if $z':=f_{\V}(z)\in \V'_{\delta y'}$. By the inductive hypothesis, for all $x\triangleright y$ there exists an isomorphism $\varphi^x:\BMP_{w}^{x}\rightarrow\!\!\!\!\!\!\!\!^{\sim}\,\,\BMP_{w'}'^{x'}$ such that $\varphi^x(s\cdot m)=f_{l,x}(s)\cdot \varphi^x(m)$, for $s\in S_k$ and $m\in \BMP_w^x$.  Moreover, if $E\not \in \E_{\delta y}$ and $x$ is on the border of $E$ with $x\triangleright y$,  by the inductive hypothesis we have an isomorphism $\varphi^E:\BMP_w^{E}\rightarrow\!\!\!\!\!\!\!\!^{\sim}\,\,\,\BMP_{w'}'^{E'}$ such that $\varphi^E(s\cdot n)=f_{l,x}(s)\cdot \varphi^E(n)$, for $s\in S_k$ and $n\in \BMP_w^E$ and such that the following diagram commutes
\begin{displaymath}\begin{xymatrix}{
\BMP_w^x \ar[r]^{\varphi^{x}}\ar[d]_{\rho_{x,E}}&\BMP_{w'}'^{x'}\ar[d]^{\rho'_{x',E'}}\\
\BMP_w^E\ar[r]_{\varphi^E}&\BMP_{w'}'^{E'}\\
}
\end{xymatrix}
\end{displaymath}

 Now, if $E:\,y\longrightarrow x $ and 
$E':\,y'\longrightarrow x'$, then
$$\BMP_w^E\cong \BMP_w^{x}/ l(E)\BMP_w^{x}\,\,\,\,\text{ and }
\,\,\,\,\BMP_{w'}'^{E'}\cong \BMP_{w'}'^{x'}/ l'(E')\BMP_{w'}'^{x'}.$$
By assumption,  $f_{l,x}(l(E))=h\cdot l'(E')$ for some invertible element $h\in k^{\times}$ and $\varphi^x(l(E)\BMP_w^x)= f_{l,x}(l(E)){\BMP'}_{w'}^{x'}=l'(E'){\BMP'}_{w'}^{x'}$. Thus the quotients are also isomorphic and so there exists  $\varphi^E:\BMP_w^E\rightarrow\!\!\!\!\!\!\!^{\sim}\,\,\BMP_{w'}'^{E'}$ such that the following diagram commutes:
\begin{displaymath}\begin{xymatrix}{
\BMP_w^x \ar[r]^{\varphi^{x}}\ar[d]_{\rho_{x,E}}&\BMP_{w'}'^{x'}\ar[d]^{\rho'_{x',E'}}\\
\BMP_w^E\ar[r]_{\varphi^E}^{\sim}&\BMP_{w'}'^{E'}\\
}
\end{xymatrix}
\end{displaymath}

Now we have to construct $\BMP_w^{\delta y}$ and $\BMP_{w'}'^{\delta y'}$. Observe that $(\varphi^x)_{x\triangleright y}$ induces an isomorphism of $S_k$-modules between the sets of sections $ \Gamma(\{  \triangleright y \},\BMP_w)\cong \Gamma(\{  \triangleright' y' \},\BMP_{w'}')$ and, from what we have observed above, the following diagram commutes:
\begin{displaymath}
\begin{xymatrix}{
\Gamma(\{  \triangleright y \},\BMP_w) \ar[d]|{\oplus_{x\triangleright y}\varphi^x }\ar[r]\ar@/^1.9pc/[rrr]<1ex>^{u_y}&\bigoplus_{x\triangleright y}\BMP_w^x \ar[d]|{ \oplus_{x\triangleright y}\varphi^x}\ar[r]&\bigoplus_{x\in\V_{\delta y}}\BMP_w^x\ar[r]^{\oplus\rho_{x,E}}\ar[d]|{\oplus_{x\in\V_{\delta y}}  \varphi_x} &\bigoplus_{E\in\E_{\delta y}}\BMP_w^E\ar[d]|{\oplus_{E\in\E_{\delta y}}\varphi^E} \\
\Gamma(\{  \triangleright' y' \},\BMP'_{w'})\ar[r] \ar@/_1.9pc/[rrr]<-1ex>_{u'_{y'}}&\bigoplus_{x' \triangleright' y'}\BMP_{w'}^{x'}\ar[r]&\bigoplus_{x'\in\V_{\delta y'}}\BMP_{w'}'^{x'}\ar[r]_{\oplus\rho'_{\!x'\!,E'}}&\bigoplus_{E\in\E_{\delta y'}}\BMP_{w'}'^{E'}
\\}
\end{xymatrix}
\end{displaymath}
 
It follows that there exists an isomorphism of $S_k$-modules $\BMP_w^{\delta y}\cong\BMP_{w'}'^{\delta y'}$  and by the unicity of the projective cover we obtain $\BMP_w^y\cong\BMP_{w'}'^{y'}$. This proves the Lemma.

\endproof

\begin{rem} Let $y,x,z,w\in \W$ be such that $y\leq w$ and $x\leq z$. If one could show that any isomorphism of posets  $[y,w]\cong [x,z]$ induces a $k$-isomorphism of moment graphs $f:\MG_{|_{[y,w]}}\rightarrow \MG_{|_{[x,z]}}$ (at least for $k=\mathbb{Q}$), then, by Lemma \ref{pullbackBMP}, the Lusztig-Dyer  \emph{Combinatorial Invariance Conjecture} (stated in \cite{Dy}) would follow. 
See \cite{BCM} for partial results on this conjecture.
\end{rem}

\subsection{Two KL-properties of the canonical sheaf}
Here we apply Lemma \ref{pullbackBMP} in order to lift some equalities concerning KL-polynomials to the moment graph setting.

From now on we denote by $\MG=(\V, \E, l, \leq)$ the Bruhat graph corresponding to a Weyl group $\W$ and $J=\emptyset$.
 As in \S\ref{bruhatMG} we denote by $\SR$ and $\T$ the set of simple reflections and all reflections, respectively, of $\W$. Recall that $\MG$ is a moment graph on the coroot lattice $\check{Q}$ and that there is a linear $\W$-action $\check{Q}$.

\subsubsection{Inverses.}  Kazhdan and Lusztig gave an inductive formula to calculate the  KL-polynomials  ((2.2.c) of \cite{KL}). From such a formula it follows easily by induction (cf. Ex.12, Chap.5 of \cite{BrBj}) that for any pair $y,w\in \W$ such that  $y\leq w$ one has  \begin{equation}P_{y,w}=P_{y^{-1},w^{-1}}.\end{equation} We translate this equality to a  $k$-isomorphism of stalks of canonical sheaves. 
\begin{lem}\label{lem_w^{-1}}Let $\W$ be a Weyl group.  The anti-involution on $\W$ defined by the mapping $x\mapsto x^{-1}$ induces a $k$-automorphism of the Bruhat graph $\MG$ for any $k$.
\end{lem}
\proof 
 The map $f_{\V}:\V\rightarrow \V$ defined by $x \mapsto x^{-1}$ is obviously a bijection. Moreover, for each $x,y\in \W$, $x\leq y$ if and only if $x^{-1}\leq y^{-1}$. So $f_{\V}:\V\rightarrow \V$ is an isomorphism of posets.

Observe that there exists a reflection $t\in \T$ such that  $y=tx$ if and only if $y^{-1}=rx^{-1}$, where $r=x^{-1}tx\in \T$. So $x-\!\!\!-\!\!\!-y\in \E$ if and only if $x^{-1}-\!\!\!-\!\!\!-y^{-1}\in \E$ . 

Thus, for every $x\in \W$ and any $v\in Y_k$, we set $f_{l,x}(v):=x^{-1}(v) $ and observe that if  $E:x-\!\!\!-\!\!\!-y=tx$, we have 
\begin{itemize}
\item[(a)] $f_{l,x}(l(x-\!\!\!-\!\!\!- tx))=x^{-1}(\check{\alpha_t})=\check{x^{-1}(\alpha_t)}=\pm l(x^{-1}-\!\!\!-\!\!\!-y^{-1})$,
 where $\pm x^{-1}(\alpha_t)\in \Delta^{\text{re}}_+$, because $x^{-1}(\alpha_t)=\pm \alpha_{x^{-1}tx}$.
\item[(b)]\begin{align*}f_{l,y}(v)\!=\!y^{-1}(v)\!=\!x^{-1}\!(tv)\!=\!x^{-1}(v)-\langle  \alpha_t, v\rangle x^{-1}(\check{\alpha_t}\!)\equiv\\
  \equiv x^{-1}(v)=f_{l,x}(v) \,\,\,\,(\text{mod}\, x^{-1}\!(\check{\alpha_t}\!))
\end{align*}
\end{itemize}
This proves that  we have a $k$-automorphism of the moment graph $\MG$ for any $k$. 
\endproof

From this we obtain the following corollary. 

\begin{cor}\label{corKLpropi} Let $w\in \W$. Then there exists an isomorphism $g:\MG_w\rightarrow \MG_{w^{-1}}$ of $k$-moment graphs on $\check{Q}$ and $\BMP_w\cong g^*\BMP_{w^{-1}}$ as $k$-sheaves on $\MG_w$ for any $k$.
\end{cor}
\proof  By Lemma \ref{lem_w^{-1}}, $f_{\V}:x\mapsto x^{-1}$ induces a $k$-isomorphism between the two complete subgraphs  $\MG_{w}$ and $\MG_{ w^{-1}}$. We may then set $g:=f_{\vert_{\MG_w}}$ and apply Lemma \ref{pullbackBMP}; the statement follows.
\endproof

\subsubsection{Multiplying by a simple reflection. Part I} Let $y,w\in \W$ and $s\in \SR$ such that $y\leq w$, $ws<w$ and $y\not\leq ws$. Under those hypotheses Kazhdan and Lusztig observed (proof of Theorem 4.2 of \cite{KL}) that 
\begin{equation}\label{prop_KL1ii}P_{y,w}=P_{ys,ws}.\end{equation}
In order to interpret (\ref{prop_KL1ii}) in our moment graph setting we need a standard combinatorial result (that actually holds for any Coxeter group):
\begin{lem}[\cite{Humph}, Lemma 7.4]\label{lem_coxeter} Let $s\in \SR$ and $v,u\in \W$ be such that $vs<v$ and $u<v$.
\begin{itemize}
\item[(i)] If $us<u$, then $us<vs$.
\item[(ii)] If $us>u$, then $us\leq v$ and $u\leq vs$.
\end{itemize}
Thus, in both cases, $us\leq v$.
\end{lem}

We are now able to define for any $k$ a  $k$-isomorphism of Bruhat (sub)graphs:

\begin{lem}\label{lem_ys} Let $y,w\in \W$ and $s\in \SR$ such that $y\leq w$, $ws<w$ and $y\not\leq ws$, then for any $k$ there is a $k$-isomorphism of moment graphs $\MG_{|_{[y,w]}}\longrightarrow\!\!\!\!\!\!\!\!^{\sim}\,\, \,\MG_{|_{[ys,ws]}}$.
\end{lem}
\proof We show that $f_{\V}:[y,w]\rightarrow[ys,ws]$, $x\mapsto xs$ is a bijection of posets inducing the identity map on the  labels.

We verify that if $x\in[y,w]$ then $xs\in[ys,ws]$. We see that $xs<x$; indeed, if it were not the case, then by Lemma \ref{lem_coxeter} (ii)  we would have $x\leq ws$, but this would imply  $y\leq ws$. In particular, this holds for $ y$,that is $ys < y$. Now, by Lemma \ref{lem_coxeter} (i); $$xs<x\text{ , }ws<w\,\Rightarrow\,\,xs\leq ws\,$$
 $$ys<y\text{ , }xs<x\,\Rightarrow\,\,ys\leq xs.$$

We now show that if $z\in[ys,ws]$ then $zs\in[y,w]$. Observe that $zs>z$; indeed, $ys<z$, $y=(ys)s>ys$ and if  $zs<z$, then by Lemma \ref{lem_coxeter} (ii), with $u=ys$ and $v=z$, we would get $y=(ys)s\leq z\leq ws$. 

Moreover, $z\leq ws<w$ and,  by Lemma \ref{lem_coxeter} (ii), $$zs>z\text{ , }ws<w\,\Rightarrow\,\,zs\leq w.$$
 $$y=(ys)s>ys\text{ , }z=(zs)s<zs\,\Rightarrow\,\,y\leq zs.$$

This completes the proof that $f_{\V}$ maps $[y,w]$ to $[ys,ws]$.

Let $x,z\in [y,w]$, then $x\leq z$ if and only if $xs\leq zs$. Indeed, we have already proved that  $xs<x$ and $zs<z$ so, by Lemma \ref{lem_coxeter} (i), with $u=x$ and $v=z$, we have $xs\leq zs$. On the other hand, $x=(xs)s>xs$ and it follows from Lemma \ref{lem_coxeter} (ii) with $u=xs$ and $v=z$ that $x=(xs)s\leq z$.

Finally from what we proved above,  for each $t\in \T$ we have that $x,tx\in[y,w]$ if and only if $xs,txs\in[ys, ws]$.This means that we have a bijection between sets of edges such that $f_{\E}(x\stackrel \gamma \rightarrow  tx )=xs\stackrel \gamma \rightarrow txs $.

Therefore $f=\left(f_{\V},(Id_{Y})_{x\in\V}\right)$ is a $k$-isomorphism of moment graphs on $\check{Q}$ for any $k$.

\endproof

So we have:

\begin{cor}\label{corKLpropii} Consider $y,w\in \W$ such that $y\leq w$ and $ws<w$, $y\not\leq ws$ for some $s\in \SR$. Let $f$ be as in Lemma \ref{lem_ys}, then $\BMP_w\cong f^*\BMP_{ws}$ as $k$-sheaves on $\MG_{|_{[y,w]}} $ for any $k$.
\end{cor}
\proof
The statement follows by combining  Lemma \ref{lem_ys} and Lemma  \ref{pullbackBMP} .
\endproof

We recollect the results of this section:

\begin{theor}\label{KLpropsthm} Let $y,w\in \W$ be such that $y\leq w$, then 
\item[(i)] $\BMP_w^y\cong \BMP_{w^{-1}}^{y^{-1}}$. 

Let $s\in \SR$ be such that $ws<w$, then 

\item[(ii)] if $y\not\leq ws$, $\BMP_w^y\cong\BMP_{ws}^{ys}$

All isomorphisms are isomorphisms of $S_k$-modules, for any $k$.
\end{theor}
\proof
\item[(i)] This follows from Corollary \ref{corKLpropi}, since two $k$-sheaves are isomorphic only if their stalks are pairwise isomorphic.
\item[(ii)] As before, the isomorphism descends from the $k$-isomorphism of  $k$-sheaves we obtained in Corollary \ref{corKLpropii}.

\endproof

\section{Invariants}\label{section_invariants}
Clearly not all equalities concerning Kazhdan-Lusztig polynomials come from $k$-isomorphisms of the underlying Bruhat graphs. In this section we develop  another technique and, as in the previous section,  we apply it in order to categorify two well-known properties of these polynomials.

\subsection{Multiplying by a simple reflection. Part II} Another property that Kazhdan and Lusztig in \cite{KL} (2.3.g)  proved is that if $y,w\in \W$ and $s\in\SR$ are such that $y\leq w$ and $ws<w$, then 
\begin{equation}\label{KLpropiii}P_{y,w}=P_{ys,w}.
\end{equation}

It is clear that in this case there is no hope of finding any $k$-isomorphism of moment graphs, since the two Bruhat intervals $[y,w]$ and $[ys,w]$ obviously have  different cardinality.

The goal of this section is to prove the following theorem.

\begin{theor}\label{Bys} For any pair $y,w\in \W$ and for any $s\in \SR$ such that $ws<w$ and $ys,y\leq w$, there exist
\begin{itemize}
\item an isomorphism of $S_k$-modules $\varphi^y:\BMP_w^y\rightarrow \BMP_w^{ys}$
\item a family of isomorphisms of $S_k$-modules $\varphi^E:\BMP_w^E\rightarrow \BMP_w^{Es}$, where $E:y-\!\!\!-\!\!\!- x\in\E$ and $Es:ys-\!\!\!-\!\!\!-xs\in\E$
\end{itemize}
 such that the following diagram commutes
\begin{equation}
\begin{xymatrix}{
\BMP_w^y\ar[r]^{\varphi^y}\ar[d]_{\rho_{y,E}}&\BMP_w^{ys}\ar[d]^{\rho_{ys,Es}}\\
\BMP_w^E\ar[r]^{\varphi^E}&\BMP_w^{Es}\\}
\end{xymatrix}
\end{equation}
and such that $\varphi^{ys}=(\varphi^y)^{-1}$.
\end{theor}

\subsection{Two preliminary lemmata}In order to prove our claim, we need two combinatorial lemmata.

Recall that 
\begin{equation*}\T=\{s_{\alpha}\,|\, \alpha\in\Delta_{+}^{\text{re}} \}=\{wsw^{-1}\,|\,w\in\W,\,s\in\mathcal{S}\}
\end{equation*} and, for all $x,y\in\W$, denote
 \begin{equation*}G_L(x,y):=\big\{  t\in \T\,|\,tx\in (x,y] \big\}
\end{equation*}
\begin{lem}\label{lem_reflns} Let $w,y\in \W$ and $s\in \SR$ be such that $y\leq w$, $ws<w$ and $ys< y$, then $$G_L(ys,w) =G_L(y,w) \cup \big\{ysy^{-1}\big\}.$$
\end{lem}
\proof We show that for all  $t\in G_ L(y,w)$ we have  $ys<tys\leq w$ as well, i.e. $t\in G_L(ys,w)$. Indeed, if $tys>ty$, then $ys<y<ty<tys$ and, by Lemma \ref{lem_coxeter} (ii) with $u=ty$ and $v=w$, $tys\leq w$ . Otherwise, $tys<ty\leq w$, $y<ty$, $ys<y$ and, by Lemma \ref{lem_coxeter} (i) with $u=y$ and $v=ty$, we obtain $ys<tys$. 

Clearly, $ysy^{-1}\in G_L(ys,w)$ and this completes the proof that the set on the righ- hand side is a subset of the one on the left.

Now we verify that if $t\in \T$, $tys\in[ys,w]$ and $ty\not \in[y,w]$, then  $t=ysy^{-1}$. Indeed, by Lemma \ref{lem_coxeter} with $u=tys$ and $v=w$, $tys\leq w$ and, if $ty\not\in[y,w]$, then $ty<y$. Moreover, $ys<y$ and so,  by Lemma \ref{lem_coxeter} (ii) with $u=ty$ and $v=y$, $tys\leq y$. So $ys<tys\leq y$ and we know that $[ys,y]=\{ys,y\}$. Thus $tys=y$, that is, $t=ysy^{-1}$.
\endproof

\begin{lem}\label{lem_interval} Let $w,y\in \W$ and $s\in \mathcal{S}$ be such that $y \leq w$, $ys<y$ and $ws<w$, then the set $[ys, w]\setminus\{ys,y\}$ is stabilized by the mapping $x\mapsto xs$.
\end{lem}
\proof Notice that $ys<y\leq w$, so it makes sense to write $[ys,w]$. Let $\mathcal{I}:=[ys, w]\setminus\{ys,y\}$ and let $x\in \mathcal{I}$. If $xs>x$, then obviously $ys<xs$ and, by  Lemma \ref{lem_coxeter} (ii) with $u=x$ and $v=w$, $xs\leq w$. On the other hand, if $xs<x$, then $xs<w$ and, by applying  Lemma \ref{lem_coxeter} (ii) with $u=ys$ and $v=x$, $ys\leq xs$. Then, in both cases $xs\in [ys, w]$ and, since $xs\neq y$ and $xs\neq ys$, we get $x\in\I$.

Finally, if $x\in \mathcal{I}$, then $xs\neq y$. Indeed $xs=y$ if and only if $x=ys\not\in \mathcal{I}$.
\endproof

\subsection{Proof of the main theorem}

Let $\ell:\W\rightarrow \mathbb{Z}_{\geq 0}$ denote the length function on $\W$. We will prove Theorem \ref{Bys} by induction on $n=\ell(w)-\ell(y)$.

If $n=0$, then $y=w$ and there is nothing to prove. If $n>0$ and $ys>y$, then $\ell(w)-\ell(ys)=n-1$ and by induction we get the desired isomorphisms.

Now, we may suppose $n>0$ and $ys<y$. Let $\I=[ys, w]\setminus \{y,ys\}$. From the inductive hypothesis, for any $x\in \I$ we get 
\begin{itemize}
\item an isomorphism of $S_k$-modules $\varphi^x:\BMP_w^x\rightarrow \BMP_w^{xs}$
\item a family of isomorphisms of $S_k$-modules $\varphi^F:\BMP_w^F\rightarrow \BMP_w^{Fs}$, where $F:x\rightarrow z\in\E^{\delta y}$ and $Es:xs\rightarrow zs\in\E^{\delta ys}$
\end{itemize}
 such that the following diagram commutes
\begin{equation}\label{diagr_comm_ys}
\begin{xymatrix}{
\BMP_w^x\ar[r]^{\varphi^x}\ar[d]_{\rho_{x,F}}&\BMP_w^{xs}\ar[d]^{\rho_{xs,Fs}}\\
\BMP_w^F\ar[r]^{\varphi^F}&\BMP_w^{Fs}\\}
\end{xymatrix}
\end{equation}
and such that $\varphi^{xs}=(\varphi^{x})^{-1}$.

Observe that our claim will follow, once we prove that there is an isomorphism of $S_k$-modules $\varphi^y:\BMP_w^y\rightarrow \BMP_w^{ys}$ compatible with the restriction maps. 
Indeed, for $E:y\rightarrow x\in\E_{\delta y}$ there exists exactly one $Es:ys\rightarrow xs\in\E_{\delta ys}$, and $\varphi^{E}$ would already have been given. If  $E:ys\rightarrow y$, then we could  set  $\varphi^E=Id$.
 Finally, for $x\neq ys$, there exists an edge $E:x\rightarrow y\in \E$ if and only if there is $Es:xs\rightarrow ys\in \E$ (cf. Lemma \ref{lem_reflns}) and in this case $\BMP_w^E\cong \BMP_w^y/l(E)\BMP_w^y\cong \BMP_w^{ys}/l(Es)\BMP_w^{ys}$, since  $E=Es$.

We will get $\varphi^y$ by defining a surjective map from $\BMP_w^y$ to $\BMP_w^{\delta ys}$.  Since $\BMP_w^{ys}$ is the projective cover of the $S_k$-module $\BMP_w^{\delta ys}$, and, since $\text{rk}_{S_k}\BMP_w^{y}\leq \text{rk}_{S_k}\BMP_w^{ys}$ (cf. Lemma 3.12. of \cite{FieNotes}), Theorem \ref{Bys} will follow from the unicity of the projective cover.

\subsubsection{Invariants} By Lemma \ref{lem_interval}, $\I$ is invariant with respect to the right multiplication by $s$ and we may define an automorphism $\sigma_s$ of the set of global sections of the Braden-MacPherson sheaf as follows. Let $m=(m_{x})\in\Gamma(\I,\BMP_w)$, then we set $\sigma_{s}(m)=(m'_x)$, where $m'_x:=\varphi^{xs}(m_{xs})$. Since the $\varphi^x$'s are, by definition, compatible with the restriction maps (see Diagram (\ref{diagr_comm_ys})), $\sigma_{s}(m)\in\Gamma(\I,\BMP_w)$. Moreover, for any $x\in \I$, $\varphi^{xs}=(\varphi^{x})^{-1}$ and so $\sigma_s$ is an involution.

Let us denote by $\Gamma^s$ the submodule of $\sigma_s$-invariant elements of $\Gamma(\mathcal{I},\BMP_w)$, and 
 by $\Gamma^{-s}$ the elements  $m\in \Gamma(\mathcal{I},\BMP_w)$ such that $\sigma_s(m)=-m$.

Let us consider $c_s:=(c_{s,x})\in \bigoplus_{x\in \W} S_k$, where $c_{s,x}:=x(\check{\alpha_s})$; then $c_s\in \Z_k$ and so it acts on  $\Gamma(\mathcal{I},\BMP_w)$ via componentwise multiplication.

\begin{lem}\label{lem_invariants} Let $(\MG_{|_{\mathcal{I}}},k)$ be a GKM-pair, then we have $\Gamma(\mathcal{I},\BMP_w)=\Gamma^s\oplus c_s \cdot\Gamma^s$.
\end{lem}
\proof (We follow \cite{Fie07a}, Lemma 2.3).

By definition, $\sigma_s $ is an involution and 2 is an invertible element in $k$, so we get $\Gamma(\mathcal{I},\BMP_w)=\Gamma^s\oplus \Gamma^{-s}$.

Let $m\in \Gamma^s$, then $\sigma_s(c_s \cdot m)=-(c_s\cdot m)$, i.e. $c_s\cdot\Gamma^s\subseteq \Gamma^{-s}$.
 Indeed,  $s(\check{\alpha_s})=-\check{\alpha_s}$ and so for any   $x\!\in\! \mathcal{I}$ we have  
\begin{equation*}
(\sigma_s(c_{s}\cdot m))_x=\varphi^{xs}(xs(\check{\alpha_s})\cdot m_{xs})=-x(\check{\alpha_s}) \cdot m_x=-c_{s,x}\cdot m_x=-(c_{s}\cdot m)_x.
\end{equation*}

We have to prove the other inclusion, that is, every element $m\in \Gamma^{-s}$ can be divided by $(x(\check{\alpha_s}))_{x\in \mathcal{I}}$ in $ \Gamma(\mathcal{I},\BMP_w)$.

If $m=(m_x)\in \Gamma^{-s}$ then $m_x=-\varphi^{xs}(m_{xs})$ and so $\rho_{xs,xs\rightarrow x}(m_{xs})=-\rho_{x,xs\rightarrow x}(m_{x})$, since the following diagram commutes:
\begin{equation}\label{diag_proof_inv}\begin{xymatrix}{
\BMP_w^{xs}\ar[r]^{\varphi^{xs}}\ar[d]|{\rho_{xs,xs\rightarrow x}}&\BMP_w^{x}\ar[d]|{\rho_{x,xs\rightarrow x}}\\
\BMP_w^{xs\rightarrow x}\ar[r]_{\varphi^{xs\rightarrow x}}&\BMP_w^{xs\rightarrow x}\\
}
\end{xymatrix}
\end{equation}

But $m$ is a section so $\rho_{xs,xs\rightarrow x}(m_{xs})=\rho_{x,xs\rightarrow x}(m_{x})$. 
It follows that $2\rho_{x,xs\rightarrow x}(m_{x})=0$, but, by definition of the canonical sheaf, 
$\ker \rho_{x,xs\rightarrow x}=\check{\alpha}_{xsx^{-1}} \BMP_w^{x}$, that is, $\check{\alpha}_{xsx^{-1}}$ divides $m_x$ in $\BMP_w^x$.

 Notice that $\check{\alpha}_{xsx^{-1}}=\pm x(\check{\alpha_s})=\pm c_{s,x}$, i.e. $c_s^{-1}\cdot m\in \bigoplus_{x\in I} \BMP_w^{x}$. We have to verify that $\rho_{x, x-\!\!\!-\!\!\!-tx}(c_{s,x}^{-1}m_{x})= \rho_{tx, x-\!\!\!-\!\!\!-tx}(c_{s,tx}^{-1}m_{tx})$ for all $t\in \T$:
\begin{eqnarray*}
&&(c_{s,tx}c_{s,x})(\rho_{tx, x-\!\!\!-\!\!\!-tx}(c_{s,tx}^{-1}m_{tx})-\rho_{x, x-\!\!\!-\!\!\!-tx}(c_{s,x}^{-1}m_x))\\
&&=c_{s,x}(\rho_{tx, x-\!\!\!-\!\!\!-tx}(m_{tx}))-c_{s,tx}(\rho_{x, x-\!\!\!-\!\!\!-tx}(m_x))\\
&&=
(c_{s,x}-c_{s,tx})\rho_{tx, x-\!\!\!-\!\!\!-tx}(m_{tx})+c_{s,tx}(\rho_{tx, x-\!\!\!-\!\!\!-tx}(m_{tx})-\rho_{x, x-\!\!\!-\!\!\!-tx}(m_x)).
\end{eqnarray*}

The term on the last line is divisible by $\check{\alpha_t}$; indeed, 
\begin{equation*}c_{s,x}-c_{s,tx}=x(\check{\alpha_s})-x(\check{\alpha_s})+\langle x(\check{\alpha_s}),\alpha_t \rangle\check{\alpha_t}\equiv 0 \,(\!\!\!\mod \alpha_t)
\end{equation*}
and 
\begin{equation*}\rho_{tx, x-\!\!\!-\!\!\!-tx}(m_{tx})-\rho_{x, x-\!\!\!-\!\!\!-tx}(m_x)=0.
 \end{equation*}

Using  the GKM-property $c_{s,tx}c_{s,x}=tx(\check{\alpha_s})\cdot x(\check{\alpha_s})$ is a multiple of $\check{\alpha_t}$ if and only if $xsx^{-1}=t$, that is $xs=tx$. Then, $m_x=-\varphi^{xs}(m_{tx})$,  $c_{s,tx}=-c_{s,x}$ and, considering that  Diagram (\ref{diag_proof_inv}) commutes, we obtain 
\begin{eqnarray*}
 \rho_{x, x-\!\!\!-\!\!\!-tx}(c_{s,t}^{-1}m_{x})&=&-c_{s,tx}^{-1} \,\rho_{x, x-\!\!\!-\!\!\!-tx}(m_{x})\\
&=&-c_{s,tx}^{-1} \,( -\rho_{tx, x-\!\!\!-\!\!\!-tx}(m_{tx}))\\
 &=& \rho_{tx, x-\!\!\!-\!\!\!-tx}(c_{s,tx}^{-1}m_{tx})
\end{eqnarray*}

Otherwise,  $xsx^{-1}\neq t$ and $\check{\alpha_t}$ divides  $\rho_{tx, x-\!\!\!-\!\!\!-tx}(c_{s,tx}^{-1}m_{tx})-\rho_{x, x-\!\!\!-\!\!\!-tx}(c_{s,x}^{-1}m_x)$ and so $\rho_{x, x-\!\!\!-\!\!\!-tx}(c_{s,x}^{-1}m_{x})= \rho_{tx, x-\!\!\!-\!\!\!-tx}(c_{s,tx}^{-1}m_{tx})$.

\endproof

\subsubsection{Building $\BMP_w^{\delta ys}$}

Let us denote 

\begin{displaymath}
\begin{xymatrix}{
\Gamma(\I,\BMP_w) \ar[r]\ar@/_1.5pc/[rrr]_{\pi_1}\ar@{^{(}->}[r]&\bigoplus_{ x\in \I}\!\BMP_w^x\ar[r]&\bigoplus_{x\in\V_{\delta y}} \BMP_w^x \ar[r]^{\oplus \rho_{x,E}}&\bigoplus_{E\in \E_{\delta y}} \BMP^E
}\end{xymatrix}
\end{displaymath}

Recall that $\BMP_w^y=u_y(\Gamma(\{>y\}, \BMP_w))$, where $u_y$ was defined as the composition of the following maps

\begin{displaymath} 
\begin{xymatrix}{
\Gamma(\{> y\},\BMP_w) \ar@{^{(}->}[r]\ar@/_1.5pc/[rrr]_{u_y}&\bigoplus_{ x>y}\!\BMP_w^x\ar[r]&\bigoplus_{x\in \V_{\delta y}}\!\!\BMP_w^x\ar[r]^{\oplus \rho_{x,E}}&\bigoplus_{E\in \E_{\delta y}} \!\!\BMP_w^E
}\end{xymatrix}
\end{displaymath}

\begin{rem}\label{rem_uy}Since $\BMP_w$ is flabby and $\I$ and $\{>y \}$ are both open sets, we get
 \begin{equation}\label{pi1=Bdy}\pi_1(\Gamma(\I,\BMP_w))=u_y(\Gamma(\{>y\}, \BMP_w))=\BMP_w^{\delta y}
\end{equation}
\end{rem}

Now, let us denote 
\begin{displaymath}
\begin{xymatrix}{
\Gamma(\I,\BMP_w) \ar[r]\ar@/_1.5pc/[rrr]_{\pi_2}\ar@{^{(}->}[r]&\bigoplus_{ x\in \I}\!\BMP_w^{x}\ar[r]&\bigoplus_{x\in\V^{\delta y}} \BMP_w^{xs} \ar[r]^{\oplus \rho_{xs,Es}}&\bigoplus_{E\in \E_{\delta y}} \BMP_w^{Es}
}\end{xymatrix}
\end{displaymath}

and define $\widetilde{\BMP_w^{\delta ys}}:=\pi_2(\Gamma(\I,\BMP_w))$.

\begin{lem}\label{lem_Bys1}\item[(i)] $\BMP_w^{\delta y}=\pi_1(\Gamma(\I, \BMP_w))={\pi_1}(\Gamma^s)$
\item[(ii)] $\widetilde{\BMP_w^{\delta ys}}=\pi_2(\Gamma(\I, \BMP_w))={\pi_2}(\Gamma^s)$
\end{lem}
\proof\item[(i)] Let $m\in\Gamma(\I,\BMP_w)$. Then, by Lemma \ref{lem_invariants}, $m=m'+c_s\cdot m''$, with $m', m''\in\Gamma^s$ and, 
if $m'=(m'_x)$, $m''=(m''_x)$,
\begin{equation*} \pi_1(m)=\big(\rho_{x,E}(m'_x)\big)_{x\in\V: y\rightarrow x\in \E}+\big(\rho_{x,E}(x(\check{\alpha_s})\cdot m''_x)\big)_{x\in\V: y\rightarrow x\in \E}
\end{equation*}

If $E:y\rightarrow x\in\E_{\delta y}$, then there exists a reflection $t\in \T$ such that $x=ty$ and we have
\begin{equation*}x(\check{\alpha_s})=ty(\check{\alpha_s})=y(\check{\alpha_s})-\langle y(\check{\alpha_s}) , \alpha_t \rangle \check{\alpha_t}
\end{equation*}

But, by definition, $\rho_{x,E}$ is a surjective map whose kernel is $l(E)\BMP_w^x=\check{\alpha_t} \BMP_w^x$ and
\begin{eqnarray*}\rho_{x,E}(x(\check{\alpha_s})\cdot m''_x)&=&\rho_{x,E}(y(\check{\alpha_s})\cdot m''_x)-
\langle y(\check{\alpha_s}) , \alpha_t \rangle \rho_{x,E}(\check{\alpha_t}\cdot m''_x)\\
&=&\rho_{x,E}(y(\check{\alpha_s})\cdot m''_x)
\end{eqnarray*}

We conclude that $\pi_1(m)=\pi_1(m'+\overline{y(\check{\alpha_s})}\cdot m'')$, where $\overline{y(\check{\alpha_s})}$ is the element of
 the structure algebra, whose components are all equal to  $y(\check{\alpha_s})$. Clearly, $m'+\overline{y(\check{\alpha_s})}\cdot m''\in\Gamma^s$ and we get the claim.

\item[(ii)] As in \textit{(i)}.
\endproof

\begin{lem}\label{lem_Bys2}There is an isomorphism of $S_k$-modules $\tau:\BMP_w^{\delta y}\rightarrow\widetilde{ \BMP_w^{\delta ys}}$ given by $(m_E)_{E\in\E_{\delta y}}\mapsto (\varphi^{E}(m_E))_{E\in\E_{\delta y}}$, that is for all $m\in\Gamma^s$, $\tau\circ\pi_1(m)=\pi_2(m)$.
\end{lem}
\proof

The element $(m_E)_{E\in\E_{\delta y}}\in\BMP_w^{\delta y}$ if and only if there exists an element $m\in\Gamma(\{>y\}, \BMP_w)$ such that $u_y(m)=(m_{E})_{E\in\E_{\delta y}}$. We have already noticed that this is the case if and only if there is an element $m'\in\Gamma(\I,\BMP)$ such that $\pi_1(m')=(m_{E})_{E\in\E_{\delta y}}$. From the previous lemma, we know that this is equivalent to the existence of an $\widetilde{m}\in\Gamma^s$ such that ${\pi_1}(\widetilde{m})=(m_E)_{E\in\E_{\delta y}}$. But, since the squares in the following diagram are all   commutative, 
\begin{displaymath}\begin{xymatrix}{
\Gamma^s\ar[d]_{Id}\ar@/^1.4pc/[rr]<1ex>^{{\pi_1}_{|_{\Gamma^s}}}\ar[r]&\bigoplus_{x\in\V_{\delta y}}\BMP_w^x\ar[r]^{\oplus \rho_{x,E}}\ar[d]|{\,\,\,\oplus \varphi^{x}}&\bigoplus_{E\in\E_{\delta y}}\BMP_w^E\ar[d]|{\,\,\,\oplus \varphi^{E}}\\
\Gamma^s\ar[r]\ar@/_1.4pc/[rr]<-1ex>_{{\pi_2}_{|_{\Gamma^s}}}&\bigoplus_{x\in\V_{\delta y}}\BMP_w^{xs}\ar[r]_{\oplus \rho_{xs,Es}}&\bigoplus_{E\in\E_{\delta y}}\BMP_w^{Es}
}
\end{xymatrix}
\end{displaymath}

we get $(\varphi^{E}(m_{E}))_{E\in\E_{\delta y}}=(\oplus \varphi^E)\circ \pi_1(\widetilde{m})=\pi_2(\widetilde{m})\in\widetilde{\BMP^{\delta ys}}$.

Analogously, $(m_{Es})_{E\in\E_{\delta y}}\in \widetilde{\BMP_w^{\delta ys}} $ if and only if $((\varphi^{E})^{-1}(m_{Es}))_{E\in\E_{\delta y}}\in\BMP_w^{\delta y}$.

\endproof

Let us denote by $\rho:\BMP_w^y\rightarrow \BMP_w^y/\check{\alpha_s}\cdot \BMP_w^y$ the canonical quotient map.

\begin{lem}We have 
\begin{equation}\BMP_w^{\delta ys}=\left\{\big(\tau\circ d_y(m_y),\rho(m_y)\big)\in\widetilde{\BMP_w^{\delta ys}}\oplus (\BMP_w^y/\check{\alpha_s}\cdot \BMP_w^y)\right\}
\end{equation}
\end{lem}
\proof
\begin{equation*}\begin{array}{rl}
\BMP_w^{\delta ys}&=u_{ys}\big(\Gamma(\{>ys\}, \BMP_w)\big)\\
&=u_{ys}\big(\left\{(m,m_y)\in\Gamma(\I, \BMP_w)\oplus \BMP_w^y\,\vert \, u_y(m_{|_{\{>y\}}})=d_y(m_y)\right\}\big)
\end{array}\end{equation*}
by Remark \ref{rem_uy}
\begin{equation*}\begin{array}{rl}
\hspace{15mm}&=u_{ys}\big(\left\{(m,m_y)\in\Gamma(\I, \BMP_w)\oplus \BMP_w^y\,\vert \, \pi_1(m)=d_y(m_y)\right\}\big)\\
&=\left\{\big(\pi_2(m),\rho(m_y)\big)\,\vert \,m\in\Gamma(\I, \BMP_w), \, m_y\in \BMP_w^y,\,\pi_1(m)=d_y(m_y)\right\}
\end{array}\end{equation*}
by Lemma \ref{lem_Bys1}
\begin{equation*}\begin{array}{rl}
\hspace{3mm}&=\left\{\big(\pi_2(m),\rho(m_y)\big)\,\vert \,m\in\Gamma^s, \, m_y\in \BMP_w^y,\,\pi_1(m)=d_y(m_y)\right\}
\end{array}\end{equation*}
by Lemma \ref{lem_Bys2}
\begin{equation*}\begin{array}{rl}
\hspace{8mm}&=\left\{\big(\!\tau\circ \pi_1(m),\rho(m_y)\big)\,\vert \,m\in\Gamma^s, \, m_y\!\in\! \BMP_w^y,\,\pi_1(m)=d_y(m_y)\right\}\\
&=\left\{\big(\!\tau\circ d_y(m_y),\rho(m_y)\big)\,\vert \, m_y\!\in\! \BMP_w^y \right\}
\end{array}\end{equation*}

\endproof

From the lemma above, it follows immediately, that there is a surjective map of $S_k$-modules $\BMP_w^y\rightarrow \BMP_w^{\delta ys}$ given by $m_y\mapsto (\tau\circ d_y(m_y), \rho(m_{y}))$ and this concludes the proof of Theorem \ref{Bys}.

\subsection{Rational smoothness and $p$-smoothness of the flag variety.} We have an easy corollary of Theorem \ref{Bys}:

\begin{cor} Let $\W$ be a finite Weyl group and $w_0$ its longest element. Then $\BMP_{w_0}^y\cong S_k$ for any $y\in \W$ and any $k$.
\end{cor}
\proof We proceed by induction on $n=\ell(w_0)-\ell(y)$. If $n=0$, by definition, $\BMP_{w_0}^{w_0}\cong S_k$. If $n\geq 1$ then there exists a simple reflection $s\in \SR$ such that $ys>y$ (so, $\ell(w_0)-\ell(ys)=n-1$). Actually, $w_0s<w_0$ for any $s\in \SR$ and, by Theorem \ref{Bys}  and inductive hypothesis, we have
$\BMP_{w_0}^{y}\cong\BMP_{w_0}^{ys}\cong S_k$.
\endproof

\begin{rem} If $k=\mathbb{Q}$ the result above corresponds to the  (rational) smoothness of flag varieties, while if $k$ is a field of characteristic $p$ it gives their $p$-smoothness (cf.\cite{FW}).  Our proof is based only on the definition of  canonical sheaf; we do not use Fiebig's multiplicity one results (see \cite{Fie06}), nor  the geometry of the corresponding flag varieties.
\end{rem}

\subsection{Parabolic setting}

Let $J\subseteq \SR$ be such that $\W_J=\langle J \rangle$ is finite with longest element $w_J$.  
Let $\W^J$ be the set of minimal representatives of the equivalence classes $\W/\W_J$. For $w\in \W^J$, denote by $\BMP_{ww_J}$, resp.$\BMP_{w}^{J}$,  the canonical sheaf on $\MG_{ ww_J}$, resp. on $\MGJ_{ w}$.  It is now easy to see that:

\begin{lem}\label{iso_latclasses}Let $\W_J$ and $w_J$ be as above and consider  $x,w\in \W^J$ such that $x\leq w$, then $\BMP_{ww_J}^x\cong\BMP_{ww_J}^{xu}$ for any $u\in W_J$.
\end{lem}
\proof We proceed by induction on  $\ell(u)$. Clearly there is nothing to prove if $\ell(u)=0$. If $\ell(u)>0$ 
then there exists an $s\in \SR$ such that $us<u$ and so by the inductive hypothesis, we get $\BMP_{ww_J}^x\cong\BMP_{ww_J}^{xus}$.
 Now for any $s\in J$, $ww_Js<ww_J$ and $xus, xu\leq ww_J$ and by Theorem \ref{Bys} we obtain the claim.
\endproof

\begin{theor}\label{BMPparab} Let $(\MG_{ww_J},k)$ be a GKM-pair and let $\W^J$ and $ w_J$ be as above. If $y,w\in \W^J$ and $y\leq w$, then there is an isomorphism of $S_k$-modules $$(\BMP_{ww_J})^{yw_J}\cong (\BMP_w^J)^y.$$
\end{theor}
\proof We proceed by induction on $n=\ell(w)-\ell(y)$. If $n=0$ the statement is trivial. Suppose we have a collection of isomorphisms of $S_k$-modules  $\eta_x :(\BMP_w^J)^x\rightarrow (\BMP_{ww_J})^{xw_J}$ for any $x$ such that $\ell(w)-\ell(x)<n$.

There is a natural injective homomorphism,
\begin{equation*}j:\Gamma(\{ >y\},\BMP_{w}^J)\rightarrow \Gamma(\{ >yw_J\},\BMP_{ww_J}), 
\end{equation*}
 defined by setting $(m_x)_{x\in(y,w]\subset \W^J}\mapsto (\widetilde{m_z})_{z\in(yw_J,ww_J]\subset \W}$, where $\widetilde{m_z}:=\psi^{z}(\eta_{x}(m_x))$ if $z\in x\W_J$ and $\psi^z:\BMP_{ww_J}^{xw_J}\rightarrow\BMP_{ww_J}^{z}$ is an isomorphism (it exists by Lemma \ref{iso_latclasses}).

We will show that such a  homomorphism induces an isomorphism  $(\BMP_{ww_J})^{\delta \,yw_J}\cong (\BMP_w^J)^{\delta y}$. Then, by the unicity of  projective cover, the statement will follow.

Let  $z\in (yw_J,ww_J]$, $z=xu$, for some $x>y\in \W^J$,  $u\in \W_J$ and $u=s_{1}\ldots s_{r}$ a reduced expression with $s_{i}\in J$ for every $i$. Moreover, let $(n_v)\in \Gamma(\{ >yw_J\},\BMP_{ww_J})$.
We prove by induction on $\ell(u)=r$ that there exists a section $(p_v)\in  \Gamma(\{ >yw_J\},\BMP_{ww_J}) $ such that $p_{xs_{1}\ldots s_i}=\psi^{xs_{1}\ldots s_i}(\eta_{x}(m_x))$ for some $m_x\in (\BMP_{w}^J)^x$ for any $i=0,\ldots, r$ and  such that $u_{yw_J}((p_v))=u_{yw_J}((n_v))$.

For the base step we have $r=0$ and there is nothing to prove.

 If $z=(xs_1 s_2\ldots s_{r-1})s_r$ then, by the inductive hypothesis, there exists a section $(q_v)\in  \Gamma(\{ >yw_J\},\BMP_{ww_J}) $ and an element $m_x\in (\BMP_{w}^J)^x$ such that $q_{xs_1\ldots s_{i}}=\psi^{xs_{1}\ldots s_i}(\eta_{x}(m_x))$ and $u_y((q_v))=u_y((n_v))$ for $i=0,\ldots, r-1$. Thus,
by Lemma \ref{lem_Bys1}, the element $(p_v)\in \bigoplus_{v>yw_J}\BMP^y$ such that $$p_{ys_1\ldots s_{r-1} s_r}=\varphi^{ys_1\ldots s_{r-1}}(p_{ys_1\ldots s_{r-1}})$$
and  
$$p_{xs_1\ldots s_{i}}=q_{xs_1\ldots s_{i}}=\psi^{xs_{1}\ldots s_i}(\eta_{x}(m_x)) \,\,\,\,\forall i<r$$
is a section on $\{>yw_J \}$ and  verifies $u_{yw_J}((\widetilde{n_v}))=u_{yw_J}((n_v))$.

Finally, from  the proof of Lemma \ref{iso_latclasses} it follows that $$\varphi^{ys_1\ldots s_{r-1}}(p_{ys_1\ldots s_{r-1}})=\varphi^{ys_1\ldots s_{r-1}}(\psi^{ys_1\ldots s_{r-1}}(\eta_x(m_{x}))=\psi^{xs_{1}\ldots s_r}(\eta_{x}(m_x)).$$

\endproof

The theorem above is the analogue of the following theorem, due to Deodhar:

\begin{theor}(\cite{Deo87})Let $\W$ be a Weyl group with $\SR$, set of simple reflections,  and $J\subseteq \mathcal{S}$ such that $\W_J$ is finite. Let $w_J$ be the longest element of $\W_J$ and $y,w\in \W^J$, then $P^{J,-1}_{y,w}=P_{yw_J,ww_J}$.
\end{theor}

\section*{Acknowledgements} I would like to thank Peter Fiebig for drawing my attention to this problem and for helpful discussions, and Rocco Chiriv\`i 
and Corrado De Concini 
for useful conversations. I would also like to thank the INdAM for funding part of my stay at the Friedrich-Alexander-Universit\"{a}t of Erlangen-N\"{u}rnberg. Finally, I am grateful to the Hausdorff Research Institute for Mathematics in Bonn for its hospitality during the program ``On the Interaction of Representation Theory with Geometry and Combinatorics", where the last version of this paper was written.

\vspace{2mm}

\end{document}